\documentclass[11pt]{amsart}
\usepackage{amsmath,amsfonts,latexsym,amssymb,amscd, graphicx,pictexwd}

\renewcommand{\P}{{\mathbb P}}
\newcommand{\E}{{\mathbb E}}
\newcommand{\R}{{\mathbb R}}

\newcommand{\PP}{{\mathbf P}}

\newcommand{\0}{{\mathbf 0}}
\newcommand{\pp}{{\mathbf p}}
\newcommand{\qq}{{\mathbf q}}
\newcommand{\cc}{{\mathbf c}}
\newcommand{\zz}{{\mathbf z}}
\newcommand{\yy}{{\mathbf y}}

\newcommand{\xx}{{\mathbf x}}
\newcommand{\nn}{{\mathbf n}}
\newcommand{\cal}{\mathcal}
\newcommand{\Co}{{\rm Co}}
\newcommand{\Cy}{{\rm Cyl}}

\newcommand{\Z}{{\mathbb Z}}
\newcommand{\calF}{{\mathcal F}}

\newtheorem{thm}{Theorem}[section]

\newtheorem{lem}[thm]{Lemma}

\numberwithin{equation}{section}

\begin{document}
\title[A shape theorem and semi-infinite geodesics for the Hammersley model]{A shape theorem and semi-infinite geodesics for the Hammersley model with random weights}

\author{Eric Cator}
\address{Delft University of Technology\\
Mekelweg 4, 2628 CD Delft, The Netherlands}
\email{E.A.Cator@tudelft.nl}

\author{Leandro P. R. Pimentel}
\address{Institute of Mathematics\\ Federal University of Rio de Janeiro}
\email{leandro@im.ufrj.br}
\thanks{Leandro P. R. Pimentel was supported by grant numbers 613.000.605 and 040.11.146 from the Netherlands Organisation for Scientific Research (NWO)}

\begin{abstract}
In this paper we will prove a shape theorem for the last passage percolation model on a two dimensional $F$-compound Poisson process, called the Hammersley model with random weights. We will also provide diffusive upper bounds for shape fluctuations. Finally we will indicate how these results can be used to prove existence and coalescence of semi-infinite geodesics in some fixed direction $\alpha$, following an approach developed by Newman and co-authors, and applied to the classical Hammersley process by W\"uthrich. These results will be crucial in the development of an upcoming paper on the relation between Busemann functions and equilibrium measures in last passage percolation models \cite{CP}.
\end{abstract}

\maketitle

\section{Introduction}\label{sec:intro}

This paper is concerned with Last Passage Percolation on a compound Poisson process, called the Hammersley process with random weights. To make this more precise, let $\PP\subseteq\R^2$ be a two-dimensional Poisson process of intensity one. On each point $\pp\in\PP$ we put a random positive weight $w_\pp$ and we assume that $\{w_\pp\,:\,\pp\in\PP\}$ is a collection of i.i.d. random variables, distributed according to a distribution function $F$, which are also independent of $\PP$. When $F$ is the Dirac distribution concentrated on $1$ (each point has weight $1$; we will denote this $F$ by $\delta_1$), then we refer to this model as the classical Hammersley model (Aldous \& Diaconis \cite{AlDi}). For each $\pp,\qq\in\R^2$, with $\pp\leq\qq$ (inequality in each coordinate), when we consider an up-right path $\varpi$ from $\pp$ to $\qq$ consisting of nondecreasing Poisson points $(\pp_1,\ldots,\pp_n)$, we will view $\varpi$ as the lowest up-right continuous path connecting all the points, starting at $\pp$ and ending at $\qq$, and then excluding $\pp$. Let $\Pi(\pp,\qq)$ denote the set of all such paths. In this probabilistic model, the last-passage time $L$ between  $\pp\leq\qq$ is defined by
$$L(\pp,\qq):=\max_{\varpi\in\Pi(\pp,\qq)}\big\{\sum_{\pp'\in\varpi\cap \PP}w_{\pp'}\big\}\,.$$
Then $L$ is super-additive,
$$L(\pp,\qq)\geq L(\pp,\zz)+L(\zz,\qq)\,.$$
A finite geodesic between $\pp$ and $\qq$ is given by the lowest path that attains the maximum in the definition of $L(\pp,\qq)$, which we will denote by $\varpi(\pp,\qq)$.

As mentioned above, the Hammersley process with random weights is a generalization of the classical Hammersley process as defined in Aldous \& Diaconis \cite{AlDi}. For this classical model, many very strong results have been proved using random matrix theory and determinantal processes, starting with the famous paper by Baik, Deift \& Johansson \cite{BaDeJo}. However, these methods rely on very specific combinatorial properties of the classical Hammersley process, which do not seem to hold in general for the Hammersley process with random weights. With this in mind, Cator and Groeneboom in \cite{CaGr1} and \cite{CaGr2} developed methods using more probabilistic arguments, in the hope that these arguments could be extended. In this paper we will prove some fundamental properties of the function $L$ and the geodesics associated to it, so that at least some of the results we have for the classical Hammersley process can indeed be extended to the general case. This will be done in several upcoming papers, starting with \cite{CP}. These fundamental properties were established for first passage percolation by Newman and co-authors (see \cite{HoNe}, \cite{LiNe}, \cite{Ne}). Their ideas were applied by W\"uthrich in \cite{Wu} to the classical Hammersley process, and we will extend these results to the Hammersley process with random weights.

In this paper we are mainly interested in two things: firstly, what is the asymptotic behavior of $L$, including its fluctuations, and secondly, can we extend the finite geodesics to semi-infinite geodesics, and can we control the fluctuations of these geodesics?

To start with the first question, we denote $\nn=(n,n)$, and define
\[F(x) = \P(w_\pp \leq x)\,\mbox{ and }\,\gamma=\gamma(F)=\sup_{n\geq 1}\frac{\E(L(\0,\nn))}{n}> 0 \,.\]
\begin{thm}\label{thm:shape}
Suppose that
\begin{equation}\label{eq:a1}
\int_0^\infty \sqrt{1-F(x)}\,dx <\infty\,.
\end{equation}
Then $\gamma(F)<\infty$ and for all $x,t>0$, as $r \to + \infty$,
\[\frac{L\left(\0,(rx,rt)\right)}{r}\to \gamma \sqrt{xt}\ \  {\rm a.s.}\ \ \ \mbox{and} \ \ \   \frac{\E L\left(\0,(rx,rt)\right)}{r} \to \gamma \sqrt{xt}\,.\]
\end{thm}

Theorem \ref{thm:shape} shows that asymptotically, $\E L$ has hyperbolic level sets, mainly due to the invariance of the Poisson process under volume preserving maps: if $\lambda >0$ and $\pp\in \R^2$, then
\begin{equation}\label{eq:sym}
\{L((x,t),(y,s))\ :\ (x,t)\leq (y,s)\} \stackrel{\cal D}{=} \left\{ L\left(\pp + (\lambda y,s/\lambda),\pp + (\lambda x,t/\lambda)\right)\right\}\,.
\end{equation}
This is because under the map $(x,t)\mapsto \pp+(\lambda x,t/\lambda)$, the distribution of the Poisson process does not change, and the up-right paths are preserved. The almost sure convergence is a standard consequence of the sub-additive ergodic theorem, once we have a bound on $\E L(\0,(r,r))$, that is linear in $r$. We will show  that (\ref{eq:a1}) is a sufficient condition to have that.

To control the fluctuations of $L$, we need more control on the distribution of the weights.
\begin{thm}\label{thm:fluctuation}
If \eqref{eq:a1} is strengthened to
\begin{equation}\label{eq:a2}
 \E e^{aw}:=\int_0^\infty \exp(a x)\,dF(x) <\infty\,\mbox{ for some }a>0,
\end{equation}
then there exist constant $c_0,c_1,c_2,c_3,c_4>0$ such that for all $r\geq c_0$
\[\P\big(|L(\0,(r,r))-\gamma r |\geq u\big)\leq c_1 \exp\Big(-c_2\frac{u}{\log(r)\sqrt{r}}\Big)\,\]
for  $u\in \left[\,c_3\sqrt{r}\log^2 r\,,\,c_4 r^{3/2}\log(r)\,\right]$.
\end{thm}
With this theorem in hand, one can actually show that, asymptotically, $L$ itself has hyperbolic level sets (shape theorem). The proof of the fluctuation result uses the method of bounded differences for martingales, following Kesten's ideas in \cite{Ke} developed for first-passage percolation models. This gives a bound on the fluctuation of $L$ around its expectation. Then adapting a clever argument used by Howard and Newman in \cite{HoNe} shows that one can replace $\E(L(\0,(r,r))$ by the shape function.

The second subject of interest to us are the geodesics. The existence of semi-infinite geodesics (or rays) for percolation like models has already been extensively study by Newman and coauthors (see Newman \cite{Ne}). They developed a general approach, based on Theorem \ref{thm:fluctuation} and on the curvature of the limit shape, that leads us to what they called the $\delta$-straightness of geodesics. This property is the key for proving the existence of rays, and will be shown in Section \ref{sec:arays}.

We need the concept of an $\alpha$-ray: for each angle $\alpha\in (0, \pi/2)$ and for each point $\xx\in\R^2$, an $\alpha$-ray starting at $\xx$ is an ordered sequence $(\pp_i)_{i\geq 0}$ in $\R^2$, with $\pp_0=\xx$, $\pp_i\in \PP$ $(i\geq 1)$ and $\pp_i\leq \pp_j$ whenever $ i\leq j$ (an up-right path). Furthermore, $\varpi(\pp_j,\pp_i)\cap \PP = \{\pp_j,\ldots,\pp_i\}$ (every  part of the path is a geodesic), and finally we must have that
$$\lim_{i\to \infty} \frac{\pp_i}{\|\pp_i\|} = (\cos\alpha, \sin\alpha)\,.$$
We will show that with probability one, every semi-infinite geodesic is an $\alpha$-ray for some $\alpha\in(0,\pi/2)$, and for every $\xx\in\R^2$ and $\alpha\in(0,\pi/2)$ there exists at least one $\alpha$-ray starting at $\xx$ (Theorem \ref{thm:existence}). Furthermore, for fixed $\alpha\in(0,\pi/2)$, with probability one, for each $\xx\in\R^2$ the $\alpha$-ray starting at $\xx$ is unique; we will denote it by $\varpi_\alpha(\xx)$. Finally, for any $\xx,\yy\in\R^2$ there exists $\cc_\alpha(\xx,\yy)$ such that $\varpi_\alpha(\xx)$ and $\varpi_\alpha(\yy)$ coalesce at $\cc_\alpha(\xx,\yy)$ (Theorem \ref{thm:unicoa}).

The proof of these results can be done by using a method introduced by Licea and Newman \cite{LiNe}, that would work in a wide context. In \cite{Wu}, W\"uthrich applied this method to the classical Hammersley model\footnote{See also Howard and Newman \cite{HoNe}, Ferrari and Pimentel \cite{FePi}} to get uniqueness and coalescence for fixed directions.

The existence, uniqueness and coalescing property of $\alpha$-rays can be used to define what is called the Busemann function: for a fixed angle $\alpha\in (0,\pi/2)$ and all $\xx,\yy\in\R^2$,
\[B_\alpha(\xx,\yy) = L(c_\alpha(\xx,\yy),\yy) - L(c_\alpha(\xx,\yy),\xx).\]
The Busemann function was also considered by Newman and co-authors for the First Passage Percolation, and by W\"uthrich for the classical Hammersley process, but more as a separate object of interest. Howard and Newman in \cite{HoNe} conjecture different scaling behavior of the Busemann function in different directions, and W\"uthrich partially answers this question in \cite{Wu}. However, in \cite{CP} we show that this Busemann function is actually closely related to equilibrium measures of a generalization of the Hammersley interacting particle system (see \cite{AlDi} and \cite{CaGr1}), called the Hammersley interacting fluid system. This leads to many interesting results for the Hammersley process with random weights. Also, in the classical Hammersley case, it allows us to give a complete specification of the scaling behavior of the Busemann function, solving the aforementioned conjecture. Furthermore, this connection is used in \cite{CP} to analyze the multi-class Hammersley process, and in an upcoming paper we will use it to determine the asymptotic speed of the second class particle in a rarefaction front. All these results rely heavily on the fundamental results in this paper, since they are essential not only for the definition of the Busemann function, but also for the necessary control of the Busemann function.\\

\paragraph{\bf Overview.} In Section \ref{sec:arays} we will introduce and prove $\delta$-straightness of semi-infinite geodesics. This is the key to obtain the existence result for $\alpha$-rays. The uniqueness and coalescence result follow from similar arguments used for the classical Hammersley case by W\"uthrich, who in turn used ideas by Licea and Newman \cite{LiNe}. In Section \ref{sec:shape} we prove Theorem \ref{thm:shape} and in Section \ref{sec:fluct} we prove Theorem \ref{thm:fluctuation}, adapting ideas from Kesten \cite{Ke} and Howard \& Newman \cite{HoNe}.

\section{Semi-infinite geodesics}\label{sec:arays}
For each $\pp\in\R^2$ and $\theta\in(0,\pi/4)$, let $\Co(\pp,\theta)$ denote the cone through the axis from $\0$ to $\pp$ and of angle $\theta$. Let $R_\0^{out}(\pp)$ be the set of points $\qq\geq \pp$ such that $\pp\in\varpi(\0,\qq)$. For fixed $\delta\in(0,1)$, we say that the geodesics starting at $\0$ in $\Co((1,1),\theta)$ are $\delta$-straight if there exist constants $M,c>0$ such that for all $\pp\in \Co((1,1),\theta)$ and $\|\pp\|\geq M$,
\begin{equation}\label{eq:delta}
R_\0^{out}(\pp)\subseteq \Co(\pp,c|\pp|^{-\delta})\,.
\end{equation}

\subsection{Controlling fluctuations of the geodesics through the curvature of the limit shape}
Our argument on how to control the fluctuations of the geodesics will very closely follow the proofs given in W\"uthrich for the classical Hammersley process. In the classical case, the control of $L(\0,(r,r))$ around its asymptotic value is stronger than our Theorem \ref{thm:fluctuation}, but our result is strong enough to extend the method to the more general Hammersley process. For details of the proof, we refer to W\"utrich.\\
For each $L>0$, let $\partial\Cy(\pp,L)$ denote the side-edge of the truncated cylinder of width $L$, that is composed of points $\qq\in\R^2$ with $\qq\geq \pp$ and $|\qq|\leq 2|\pp|$, and such that the Euclidean distance between $\qq$ and the line through $\0$ and $\pp$ equals $L$. Assume that $\pp=(x,t)\in \Co((1,1),\theta)$, that $\qq\in R^{out}$ and that $\qq\in \partial\Cy(\pp,|\pp|^{1-\delta})$. Then
$$L(\0,\qq)=L(\0,\pp)+L(\pp,\qq)\,$$
or, equivalently,
$$f(\qq)-f(\qq-\pp)-f(\pp)=\Delta(\0,\pp)+\Delta(\pp,\qq)-\Delta(\0,\qq)\,,$$
where $f(\pp)=f(x,t)=\gamma\sqrt{xt}$ is the shape function and
$$\Delta(\pp,\qq)=L(\pp,\qq)-f(\qq-\pp)\,.$$
 On the other hand, since $\qq\in \partial\Cy(\pp,|\pp|^{1-\delta})$,
$$f(\qq)-f(\qq-\pp)-f(\pp)\geq c_0|\pp|^{1-2\delta}\,,$$
for a finite constant $c_0>0$, depending on $\theta\in (0,\pi/4)$ (here we use the curvature of the shape function; see Lemma 2.1 in W\"uthrich). Notice that if $\delta\in(0,1/4)$ then $1-2\delta\in(1/2,1)$ and so
$$|\pp|^{1-2\delta}>> |\pp|^{1/2}\log|\pp|\,.$$
Hence, by Theorem \ref{thm:fluctuation}, for $\delta\in(0,1/4)$, we must have that if $\qq\in \partial\Cy(\pp,|\pp|^{1-\delta})$, then with very high probability $\qq\not\in R^{out}$. This can be formalized to prove the following lemma:
\begin{lem}\label{delta1}
Fix $\delta\in(0,1/4)$ and $\theta\in(0,\pi/4)$. For each $\pp=(x,t)\in \Co((1,1),\theta)$ and $\qq\in \partial \Cy(\pp,|\pp|^{1-\delta})$, let $G_\delta(\pp,\qq)$ be the event that $\pp\in\varpi(\0,\qq)$. Then there exist finite constants $\kappa,c_1,c_2>0$ such that
$$\P\left(G_\delta(\pp,\qq)\right)\leq c_1e^{-c_2|\pp|^{\kappa}}\,.$$
\end{lem}
We extend this Lemma to hold uniformly for $\qq$ and $\pp$ in a fixed-size finite box, then use the boxes around $\qq$ to cover the side-edge of the cylinder to get:
\begin{lem}\label{delta2}
Fix $\delta\in(0,1/4)$ and $\theta\in(0,\pi/4)$. For each $\pp=(x,t)\in \Co((1,1),\theta)$, let $G_\delta(\pp)$ be the event that there exists $\pp'\in \pp+[0,1]^2$ and $\qq\in \partial \Cy(\pp,|\pp|^{1-\delta})$ such that $\pp'\in\varpi(\0,\qq)$. Then there exist finite constants $\kappa,c_3,c_4>0$ such that
$$\P\left(G_\delta(\pp)\right)\leq c_3e^{-c_4|\pp|^{\kappa}}\,.$$
\end{lem}
Now we can show $\delta$-straightness by ``gluing'' together these cylinders: if a geodesic starts at close to $\pp$, with high probability it will exit the bottom edge of the cylinder $\Cy(\pp,|\pp|^{1-\delta})$. Then we cover this bottom edge with boxes $\pp_2+[0,1]^2$, where $|\pp_2|\geq 2|\pp|$, and for each of these $\pp_2$ we consider the cylinder $\Cy(\pp_2,|\pp_2|^{1-\delta})$, and so on. With Borel-Cantelli we can make the probability that a geodesic starting close to $\pp$ will ever leave through the outer edges of the boundary cylinders very small. The cylinders at the next step of the procedure have a slightly different angle than the cylinders in the previous step, but the changes in these angles are bounded by a geometric series, which means that all cylinders are contained in a cone starting at $0$, of angle $|\pp|^{-\delta}$. This is basically the same argument used for the proof of Lemma 2.4 in\cite{Wu} in the classical set-up.  This leads us to
\begin{lem}\label{delta3}
Fix $\delta\in(0,1/4)$ and $\theta\in(0,\pi/4)$. There exist $c_0,c_1,\kappa,M>0$ such that for all $\pp\in \Co((1,1),\theta)$ with $|\pp|>M$, we have
\[\P\left(\left(\bigcup_{\pp'\in \pp+[0,1]^2} R_\0^{out}(\pp')\right)\subset \Co(\pp,|\pp|^{-\delta})\right)\geq 1-c_0e^{-c_1|\pp|^\kappa}\,.\]
Furthermore, with probability one, there exists $M>0$ such that for all $\pp\in \Co((1,1),\theta)$ with $|\pp|\geq M$,
\[R_\0^{out}(\pp)\subset \Co(\pp,|\pp|^{-\delta}).\]
\end{lem}

\subsection{Existence, uniqueness and coalescence of $\alpha$-rays}

With Lemma \ref{delta3} in hands, one can show existence of $\alpha$-rays. The proof of the next theorem follows mutatis mutandis the proof of Theorem 3.4 of \cite{Wu} (compare Lemma 2.4 of \cite{Wu} with our Lemma \ref{delta3}).
\begin{thm}\label{thm:existence}
With probability one, every ray  is an $\alpha$-ray for some $\alpha\in(0,\pi/2)$, and for every $\pp\in\PP$ and $\alpha\in(0,\pi/2)$ there exists at least one $\alpha$-ray starting at $\pp$.
\end{thm}

Uniqueness and coalescence of $\alpha$-rays do not depend upon $\delta$-straightness. The proof of these can be done by using a method  introduced by Licea and Newman \cite{LiNe}, that would work in a wide context. In \cite{Wu}, W\"uthrich applied this method to the classical Hammersley model\footnote{See also Howard and Newman, Ferrari and Pimentel} to get uniqueness and coalescence for fixed directions. Here we state without proof the analogous result for the Hammersley model with random weights. The reader can convince her- or himself of the validity of the theorem by checking that the proof given by W\"uthrich can be adapted mutatis mutandis to our set-up.
\begin{thm}\label{thm:unicoa}
For fixed $\alpha\in(0,\pi/2)$, with probability one, for each $\pp\in\R^2$ there exists a unique $\alpha$-ray starting at $\pp$, which we denote by $\varpi_\alpha(\pp)$. Furthermore, for any $\pp,\qq\in\R^2$ there exists $\cc_\alpha(\pp,\qq)$ such that $\varpi_\alpha(\pp)$ and $\varpi_\alpha(\qq)$ coalesce at $\cc_\alpha(\pp,\qq)$.
\end{thm}

\section{Proof of Theorem \ref{thm:shape}}\label{sec:shape}
Equation \eqref{eq:sym} shows that it is enough to prove the theorem for $(x,t)=(1,1)$ (choose $\lambda = \sqrt{t/x}$). When considering only one ray, the convergence of $L$ is a standard consequence of Liggett's version of the superadditive ergodic theorem \cite{Li}, as soon as we can show that
\[ \limsup_{r \to \infty} \frac{\E L\left(\0,(r,r)\right)}{r} <\infty\,.\]
For each $p\in[0,1]$ denote by $\E_p$ expectation for the Hammersley last passage model induced by Bernoulli weights $w'_{\pp}$, where $\P(w'_\pp = 1) = p$. This coincides with the classical Hammersley model, but with Poisson intensity $p$ (instead of $1$). From \eqref{eq:sym}, it is well known (see Aldous \& Diaconis (1995) or Cator \& Groeneboom (2005)) that
\[ \lim_{r \to \infty} \frac{\E_pL\left(\0,(r,r)\right)}{r} = \gamma(1)\sqrt{p}\,,\]
for some $\gamma(1)<\infty$ (our notation). Now we use an idea introduced in Martin (2004):
\begin{eqnarray*}
L(\0,\pp) & = & \max_{\varpi\in\Pi(\0,\pp)}\big\{\sum_{\pp'\in\varpi}w_{\pp'}\big\}\\
& = & \max_{\varpi\in\Pi(\0,\pp)}\big\{ \int_0^\infty \sum_{\pp'\in\varpi}1_{\{w_{\pp'}>x\}} \,dx\big\}\\
& \leq & \int_0^\infty \max_{\varpi\in\Pi(\0,\pp)}\big\{\sum_{\pp'\in\varpi}1_{\{w_{\pp'}>x\}} \big\}\,dx\,.
\end{eqnarray*}
The integrand in the last line corresponds to the Bernoulli model with $p=1-F(x)$. This means that
\begin{eqnarray*}
\limsup_{r \to \infty} \frac{\E L\left(\0,(r,r)\right)}{r} & \leq & \limsup_{r \to \infty} \int_0^\infty \frac{\E_{1-F(x)}L\left(\0,(r,r)\right)}{r}\, dx \\
& = & \gamma(1) \int_0^\infty \sqrt{1-F(x)}\,dx\,.
\end{eqnarray*}

\section{Proof of Theorem \ref{thm:fluctuation}}\label{sec:fluct}

The proof of Theorem \ref{thm:fluctuation} follows Kesten's approach developed for first-passage times in lattice firs-passage percolation models \cite{Ke}. It is based on the method of bounded increments applied to $L$.
\begin{lem}\label{Kesten}
Let $\{\calF_k\}_{0\leq k\leq N}$ be a filtration and let $\{U_k\}_{0\leq k\leq N}$ be a family of positive random variables that are $\calF_N$ measurable. Let $\{M_k\}_{0\leq k\leq N}$ be a martingale with respect to $\{\calF_k\}_{0\leq k\leq N}$. Assume that the increments $\Delta_k:=M_k-M_{k-1}$ satisfy
\begin{equation}\label{hip1}
|\Delta_k|\leq c\mbox { for some }c>0
\end{equation}
and
\begin{equation}\label{hip2}
\E \left(\Delta_k^2\mid\calF_{k-1}\right)\leq \E \left(U_k\mid\calF_{k-1}\right)\,.
\end{equation}
Assume further that for some constants $0<c_1,c_2<\infty$ and $x_0\geq e^2c^2$ we have
\begin{equation}\label{hip3}
\P\left(\sum_{k=1}^N U_k>x\right)\leq c_1 \exp(-c_2 x)\mbox{ when }x\geq x_0\,.
\end{equation}
Then irrespective of the value of $N$, there exists universal constants $0<c_3,c_4<\infty$ that do not depend on $N,c,c_1,c_2$ and $x_0$, nor on the distribution of $\{M_k\}_{0\leq k\leq N}$ and $\{U_k\}_{0\leq k\leq N}$, such that
\begin{equation}\label{bound1}
\P(M_N-M_0\geq x)\leq c_3\left\{1+c_1+\frac{c_1}{c_2x_0}\right\}\exp\left(-c_4\frac{x}{\sqrt{x_0}}\right)\,
\end{equation}
whenever $x\leq c_2 \,x_0^{3/2}$.
\end{lem}
\noindent{\bf Proof:} See Theorem 3 in \cite{Ke}. \hfill $\Box$\\

We decompose $L$ as a sum of martingales increments as follows. For each integer $r\geq 1$, let $L_r:=L(\0,(r,r))$ and consider a partition of the two dimensional square $[0,r]^2=\cup_{l=1}^{N}B_l$ into $N=r^2$ disjoint squares of size one. Let $\calF_0=\{\emptyset, \Omega\}$ and for each $k=1,\dots,N$ consider the $\sigma$-algebra
$$\calF_k=\sigma\left(\{\,\omega_\pp\,:\,\pp\in\cup_{l=1}^{k}B_l\cap\PP\}\right)\,,$$
and the Doob martingale
$$M_k:=\E\left(L_r\mid\calF_k\right)\,.$$
Denote by $\P_l$ the probability law induced by  $\{\,\omega_\pp\,:\,\pp\in B_l\cap\PP\}$, and by $\Omega_l$ the underlying sample space. For $\omega,\sigma\in\prod_{l=1}^{N}\Omega_l$ let $[\omega,\sigma]_k:=(\omega_1,\dots,\omega_k,\sigma_{k+1},\dots,\sigma_{r^2})\in\prod_{l=1}^{N}\Omega_l$. Then
\begin{eqnarray}
\nonumber\Delta_k(\omega_1,\dots,\omega_k)&:=&M_k-M_{k-1}\\
\nonumber&=&\int L_r[\omega,\sigma]_{k}\prod_{l=k+1}^{r^2}d\P_l(\sigma_l)-\int L_r[\omega,\sigma]_{k-1}\prod_{l=k}^{N}d\P_l(\sigma_l)\\
\nonumber&=&\int L_r[\omega,\sigma]_{k}-L_r[\omega,\sigma]_{k-1}\prod_{l=k}^{N}d\P_l(\sigma_l)\,
\end{eqnarray}
(To integrate $M_k$ over $\sigma_k$ does not change it, and it allows us to put $M_k$ and $M_{k-1}$ under the same integral.).
For each $1\leq k\leq N$ let
\begin{equation}
\nonumber Z_k:=\sum_{\pp\in B_k\cap\PP}w_\pp \,.
\end{equation}
Then $\{Z_k\}_{1\leq k\leq N}$ is an i.i.d. collection of random variables such that
$$\E e^{aZ_1}=\exp\left(\E e^{aw}-1\right)<\infty\,.$$
Let $I_k$ denote the indicator function of the event that the geodesic $\varpi_r:=\varpi(\0,(r,r))$ has a point $\pp\in B_k\cap\PP$.
\begin{lem}
\begin{equation}\label{increm}
|L_r[\omega,\sigma]_{k}-L_r[\omega,\sigma]_{k-1}|\,\leq\,\max\left\{I_k[\omega,\sigma]_k,I_k[\omega,\sigma]_{k-1}\right\}\times\max\left\{ Z_k(\omega_k),Z_k(\sigma_k)\right\}\,.
 \end{equation}
 \end{lem}
\noindent{\bf Proof:}
We note that there will be no difference between $L_r[\omega,\sigma]_{k}$ and $L_r[\omega,\sigma]_{k-1}$, if no geodesic has a point $\pp\in B_k\cap\PP$ (recall that $[\omega,\sigma]_{k}$ and $[\omega,\sigma]_{k-1}$ only differ inside $B_k$), and this corresponds to the first factor in the right hand side of  \eqref{increm}. And, if one of them does intersect, then the increment can not be greater then the total weight inside the box $B_k$. (Compare with (2.12) in \cite{Ke}.)
\hfill $\Box$\\

The next step is to construct, from $\{w_\pp:\pp\in\PP\}$, a new process $\{\bar{w}_\pp:\pp\in\bar\PP\}$, by truncating the original model inside each box $B_k$, if $Z_k> b\log r$, in order to have
$$\sum_{\pp\in B_k\cap\bar\PP}\bar{w}_\pp \leq b\log r\,.$$
We do this truncating by multiplying all weights in box $B_k$ with an appropriate (small) factor. Note that the truncated process in each box is still independent of all the other boxes. Let us denote the configuration induced by the truncated model by $\bar{\omega}$ and for each random variable $Y$ that depends on $\omega$, let us write $\bar{Y}(\omega):=Y(\bar\omega)$. Thus,
$$|\bar L_r[\omega,\sigma]_{k}-\bar L_r[\omega,\sigma]_{k-1}|\leq (b\log r)\max\left\{\bar I_k[\omega,\sigma]_k,\bar I_k[\omega,\sigma]_{k-1}\right\}\,,$$
and hence
\begin{eqnarray}
\label{truncincrem1} |\bar\Delta_{k}(\omega_1,\dots,\omega_{k})|&\leq& b\log r\int\max\left\{\bar I_k[\omega,\sigma]_k,\bar I_k[\omega,\sigma]_{k-1}\right\} \prod_{l=k}^{N}d\P_l(\sigma_l)\\
\label{truncincrem2} &\leq& b\log r\,.
\end{eqnarray}
The upper bounds  \eqref{truncincrem1} and \eqref{truncincrem2} allow us to apply Lemma \ref{Kesten} to get concentration inequalities for $\bar M_{N}-\bar M_0=\bar L_r(\omega)-\E \bar L_r(\omega)$.
\begin{lem}\label{conditions}
Let $U_k:=2(b\log r)^2 I_k$. Then $\bar\Delta_k\leq b\log r$ and
$$\E(\bar\Delta_k^2\mid\calF_{k-1})\leq \E(\bar U_k\mid\calF_{k-1})\,.$$
\end{lem}
\noindent{\bf Proof:}
\begin{eqnarray*}
\E(\bar\Delta_k^2\mid\calF_{k-1})&&=\int\left\{\int \bar L_r[\omega,\sigma]_{k} - \bar L_r[\omega,\sigma]_{k-1}\prod_{l=k}^{N}d\P_l(\sigma_l)\right\}^2d\P_k(\omega_k)\\
\leq&&\int\left\{\int \max\left\{\bar I_k[\omega,\sigma]_k,\bar I_k[\omega,\sigma]_{k-1}\right\}\times\max\left\{ \bar Z_k(\omega_k),\bar Z_k(\sigma_k)\right\}\prod_{l=k}^{N}d\P_l(\sigma_l)\right\}^2d\P_k(\omega_k)\\
\leq&&\int\int \max\left\{\bar I_k[\omega,\sigma]_k,\bar I_k[\omega,\sigma]_{k-1}\right\}\times\max\left\{ \bar Z_k(\omega_k),\bar Z_k(\sigma_k)\right\}^2\prod_{l=k}^{N}d\P_l(\sigma_l)d\P_k(\omega_k)\\
\leq&&\int\int \left(\bar I_k[\omega,\sigma]_k + \bar I_k[\omega,\sigma]_{k-1}\right)(b\log r)^2\prod_{l=k}^{N}d\P_l(\sigma_l)d\P_k(\omega_k)\\
= && \E(\bar U_k\mid\calF_{k-1})\,.
\end{eqnarray*}
 \hfill $\Box$\\

\begin{lem}\label{optimalpath}
Let $|\varpi(\pp,\qq)|:=\#(\varpi(\pp,\qq)\cap\PP)$. For each integer $r\geq 1$ and $x\geq2(\log 2+ 2e)r$ we have that
$$\P\left(|\varpi_r(\bar\omega)|>x\right)\leq  \exp\left(-\frac{x}{2}\right)\,.$$
\end{lem}
\noindent The relevance of this Lemma is of course that, using $|\varpi_r|=|\bar \varpi_r|$, we conclude
\begin{equation}\label{eq:sumUk}
\sum_{k=1}^N\bar U_k \leq (b\log r)^2|\varpi_r|.
\end{equation}
\noindent{\bf Proof of Lemma \ref{optimalpath}:}
For each $i,j\in\{1,\dots,r\}$ and $\zz=(i,j)\in\Z^2$, let
$$X_\zz:=\#(B_\zz\cap\PP)\,,$$
where $B_\zz:=\zz+[0,1]^2$. Let $\Gamma_r$ be the set of all up-right $\Z^2$ lattice paths from $(0,0)$ to $(r-1,r-1)$ and
$$G_r=G_r(\omega):=\max_{\gamma\in\Gamma_r}\sum_{\zz\in\gamma}X_\zz\,.$$
Then $|\varpi_r(\bar\omega)|\leq G_r(\omega)$ and hence
\begin{equation}\label{latticepath}
\P\left(|\varpi_r(\bar\omega)|>x\right)\leq \P\left(G_r(\omega) > x\right)\,.
\end{equation}
Now, $|\Gamma_r|=2^r$ and for each fixed path $\gamma\in\Gamma_r$,
$$\P\left(\sum_{\zz\in\gamma}X_\zz>x\right)=\P\left(\sum_{l=1}^{2r}X_l>x\right)\,,$$
where $X_l$ for $l=1,\dots,r$ are i.i.d. Poisson random variables of intensity $1$. Thus, by Markov's inequality,
\begin{eqnarray}
\nonumber\P\left(G_r(\omega) > x\right)&\leq&2^r\P\left( \sum_{l=1}^rX_l> x\right)\\
\nonumber &\leq&2^r e^{-x}(\E e^{X_1})^{2r}\\
\nonumber&=&\exp\left(-x+(\log 2+2\log\E e^{X_1})r\right)\\
\nonumber&=&\exp\left(-x+(\log 2+2e)r\right)\\
\nonumber&\leq& \exp\left(-\frac{x}{2}\right)\,,
\end{eqnarray}
if $x\geq 2(\log 2+ 2e)r$.
 \hfill $\Box$\\

By choosing $b>0$ large enough, one shows that the truncated model is a good approximation of the original model, in the sense that the probability that they will differ by $u$ goes exponentially fast to zero in $u$ (Compare with (2.30) and (2.34) in \cite{Ke}).
\begin{lem}\label{lem:trunc}
Let $b=6/a$ and $r\geq \E e^{a Z_1}/\log 2$. Then
$$\P\left(L_r(\omega)-\bar L_r(\omega)>x\right)\leq 2 \exp\left(-\frac{a}{2}x\right)\,.$$
\end{lem}

\noindent{\bf Proof:} Fix $b>0$ and a positive integer $r\geq 1$ ($N=r^2$). Notice that
\begin{equation}\label{comparing}
0\leq L_r(\omega)-\bar L_r(\omega)\leq \sum_{l=1}^{N} Z_lI\left\{Z_l>b\log r\right\}\,.
\end{equation}
By Markov's inequality,
\begin{eqnarray}
\nonumber \P\left(\sum_{l=1}^{N}Z_l I\{Z_l> b\log r\} > x\right)&\leq& e^{-\frac{a}{2}x}\Big[\E\left(e^{\frac{a}{2}Z_1I\{\frac{a}{2}Z_1>\frac{a b}{2}\log r\}}\right)\Big]^{N}\\
\nonumber&=&  e^{-\frac{a}{2}x}\Big[\E\left(e^{\frac{a}{2}Z_1I\{e^{a Z_1}>r^{\frac{a b}{2}}e^{\frac{a Z_1}{2}}\}}\right)\Big]^{N}\,.
\end{eqnarray}
On the other hand,
$$e^{\frac{a}{2}Z_1I\{e^{aZ_1}>r^{\frac{ab}{2}}e^{\frac{aZ_1}{2}}\}}\leq 1+ e^{\frac{a}{2}Z_1}I\{e^{aZ_1}>r^{\frac{ab}{2}}e^{\frac{aZ_1}{2}}\}\leq 1 +\frac{e^{aZ_1}}{r^{\frac{ab}{2}}}\,,$$
and hence,
$$ \P\left(\sum_{l=1}^{N}Z_l I\{Z_l> b\log r\} > x\right)\leq e^{-\frac{a}{2}x}\Big[1+\frac{\E e^{aZ_1}}{r^{\frac{ab}{2}}}\Big]^{N}\,.$$
Now,
\begin{eqnarray}
\nonumber \log\left(\left[1+\frac{\E e^{aZ_1}}{r^{\frac{ab}{2}}}\right]^{N}\right)&=&r^2\log\left(1+\frac{\E e^{aZ_1}}{r^{\frac{ab}{2}}}\right)\\
\nonumber&\leq&r^2\frac{\E e^{aZ_1}}{r^{\frac{ab}{2}}}\\
\nonumber&=&r^{\frac{4-ab}{2}}\E e^{aZ_1}\\
\nonumber&\leq&\log 2\,.
\end{eqnarray}
if we take $b=6/a$ and $r\geq \E e^{a Z_1}/\log 2$. Together with \eqref{comparing}, this proves Lemma \ref{lem:trunc}. \hfill $\Box$\\

\begin{lem}\label{concentration}
If we assume \eqref{eq:a2}, then there exist constants $b_0,b_1,b_2,b_3>0$ such that for all $r\geq b_0$
\begin{equation}\label{e-talagrand}
\P\left(|L\left(\0,(r,r)\right)-\E L\left(\0,(r,r)\right)|\geq u\right)\leq b_1\exp\left(-b_2\frac{u}{\log(r)\sqrt{r}}\right)\,,
\end{equation}
 for $u\in(0,b_3r^{3/2}\log(r)]$.
\end{lem}
\noindent{\bf Proof:} We have checked all the conditions of Lemma \ref{Kesten} applied to the truncated process, where we take for a large enough constant $C>0$, $x_0 = C(b\log r)^2r$, $c=b\log(r)$, $c_1=1$ and $c_2=1/(2b\log r)^2$ (this follows from Lemma \ref{conditions}, Lemma \ref{optimalpath} and Equation \eqref{eq:sumUk}). So there exist $c_3, c_4>0$ such that
\[ \P\left(\bar M_N - \bar M_0 \geq u\right) \leq c_3\left\{2+\frac{4}{Cr}\right\}\exp\left(-c_4\frac{u}{(b\log r)\sqrt{Cr}}\right),\]
for all $u\leq \frac14C^{3/2}(b\log r)r^{3/2}$. Using Lemma \ref{lem:trunc} and the fact that $L_r\geq \bar L_r$,we can see that there exists $M>0$, such that for all $r\geq 1$, $|\E(L_r) - \E(\bar L_r)|\leq M$. Therefore, for $u\geq M$,
\begin{eqnarray*}
\P(|L_r - \E(L_r)|\geq 2u) & \leq & \P(|L_r - \E(\bar L_r)|\geq 2u - M) \\
& \leq & \P(|\bar L_r - \E(\bar L_r)|\geq u) + \P(|L_r - \bar L_r|\geq u - M)\\
& = & \P\left(\bar M_N - \bar M_0 \geq u\right) + \P(L_r - \bar L_r\geq u - M).
\end{eqnarray*}
Again using Lemma \ref{lem:trunc}, we can choose $b_0 = \E e^{a Z_1}/\log 2$, $b_1>0$ large enough and $b_2,b_3>0$ small enough such that \eqref{e-talagrand} holds not only for $u\in [2M,b_3r^{3/2}\log(r)]$, but also for $0\leq u\leq 2M$.
\hfill $\Box$\\

\begin{lem}\label{HoNe}
There exists a constant $c_0>0$ such that
$$\gamma r - c_0 r^{1/2}\log^2 r\leq \E L\left(\0,(r,r)\right) \leq \gamma r\,.$$
\end{lem}
\noindent{\bf Proof:} We note that the right hand side of \eqref{HoNe} follows from the definition of $\gamma$. To prove that the left hand side of \eqref{HoNe} also holds, we parallel Howard \& Newmann \cite{HoNe}. We start by noting that it is enough to prove \eqref{HoNe} for integer values of $r$. Denote by $H_r$ the set of points $(x,t)$ such that $x,t\geq 0$ and $x+t=r$. Then
$$L(\0,(2r,2r))\leq \max_{\xx\in H_{2r}}L(\0,\xx) + \max_{\xx\in H_{2r}}L(\xx,(r,r))\,,$$
and hence, by symmetry with respect to $H_r$,
$$\E L\left(\0,(2r,2r)\right)\leq 2\E \max_{\xx\in H_{2r}} L(\0,\xx)\,.$$
Define, for $k\in \{0,1,\ldots, 2r\}$, $\xx_k=(k,2r-k)$, and for $k\in \{1,\ldots,2r\}$, $\zz_k=(k-1,2r-k)$. For any $\xx\in H_r$, there exists a $k\in \{1,\ldots,2r\}$ such that $\zz_k\leq \xx$. We have
\[ L\left(0,\xx\right) \leq \max\left\{L(\0,\xx_{k-1})\,,\,L(\0,\xx_k)\right) + L(\zz_k,\zz_k+(1,1)\}.\]
This implies that
\[ \E L\left(\0,(2r,2r)\right)\leq 2\E\max_{0\leq k\leq 2r}L(\0,\xx_k) + 2\E\max_{1\leq k\leq 2r}L(\zz_k,\zz_k+(1,1)).\]
The second term on the righthand side is bounded by the expectation of the maximum of the total weight in the $2r$ squares $[\zz_k,\zz_k+(1,1)]$, which is clearly bounded by $c\log(r)$, for some constant $c>0$ (here we use \eqref{eq:a2}). So we get
\[ \E L\left(\0,(2r,2r)\right)\leq 2\E\max_{0\leq k\leq 2r}L(\0,\xx_k) + c\log(r).\]
By \eqref{eq:sym},
$$\max_{\xx\in H_{2r}}\E L(\0,\xx)=\max_{x\in[0,2r]}\E L\left(\0,(\sqrt{2r-x}\sqrt{x},\sqrt{2r-x}\sqrt{x})\right)=\E L\left(\0,(r,r)\right)\,,$$
and thus
\begin{eqnarray}
\nonumber\E L\left(\0,(2r,2r)\right)&\leq& 2\E \max_{0\leq k\leq 2r} L(\0,\xx_k) + c\log(r)\\
\nonumber&\leq& 2\E \max_{0\leq k\leq 2r} \left\{L(\0,\xx_k)-\E L(\0,\xx_k)\right\} +2\max_{0\leq k\leq 2r} \E L(\0,\xx_k) + c\log(r)\\
\label{eq:sub1}&\leq& 2\E \max_{0\leq k\leq 2r} \left\{L(\0,\xx_k)-\E L(\0,\xx_k)\right\} +2\E L(\0,(r,r)) + c\log(r)\,.
\end{eqnarray}
Now define
\[ M_r = \max_{0\leq k\leq 2r} \left\{L(\0,\xx_k)-\E L(\0,\xx_k)\right\}.\]
Define for a large constant $C>0$, the event
\[ A = \{ L(\0,\xx_k)-\E L(\0,\xx_k)\leq C\log^2(r)\sqrt{r}\ (\forall\ 0\leq k\leq 2r)\}.\]
Then
\[ M_r \leq C\log^2(r)\sqrt{r}1_A + L(\0,(2r,2r))1_{A^c}.\]
Therefore,
\[ \E M_r \leq C\log^2(r)\sqrt{r} + \sqrt{\E \left[L(\0,(2r,2r))^2\right]\cdot \P(A^c)}.\]
We crudely bound $L(\0,(2r,2r))$ by the total weight in the square $[\0,(2r,2r)]$, and see that there exists a constant $c_1>0$ such that
\[ \E \left[L(\0,(2r,2r))^2\right] \leq c_1^2r^4.\]
We can use Lemma \ref{concentration} to conclude that
\begin{eqnarray*}
\P(A^c) & \leq & \sum_{k=0}^{2r}\P\left(L(\0,\xx_k)-\E L(\0,\xx_k)> C\log^2(r)\sqrt{r}\right)\\
& \leq & (2r+1)b_1\exp\left(-b_2C\log(r)\right).
\end{eqnarray*}
By increasing $C$, this shows that there exists $c_2>0$ such that for all $r\geq 1$,
\[ \E M_r\leq c_2\log^2(r)\sqrt{r}.\]
Together with \eqref{eq:sub1}, this proves that there exists $b>0$ such that for all $r\geq 1$
\begin{equation}\label{eq:sub2}
\E L\left(\0,(2r,2r)\right)-b r^{1/2}\log^2 r\leq 2\E L\left(\0,(r,r)\right)\,.
\end{equation}
By Lemma 4.2 of \cite{HoNe}, \eqref{eq:sub2} implies that the left hand side of \eqref{HoNe} is true.

\hfill $\Box$\\
The results of Lemma \ref{concentration} and Lemma \ref{HoNe} now easily combine to Theorem \ref{thm:fluctuation}.

\end{document}